\documentclass{article}
\usepackage[utf8]{inputenc}
\usepackage{amsmath,amsthm,amscd,amssymb,latexsym}

\usepackage[english]{babel}
\usepackage{setspace}
\usepackage[numbers]{natbib}
\usepackage[svgnames]{xcolor}
\usepackage[linkcolor=Blue,colorlinks=true,citecolor=ForestGreen,filecolor=red,urlcolor=Blue]{hyperref}
\usepackage{cleveref}[2011/12/24]
\usepackage[pagewise,displaymath,mathlines]{lineno}
\usepackage{scalerel}
\usepackage{fullpage}

\usepackage{tikz}

\usepackage{authblk}

\newtheorem{theorem}{Theorem}
\newtheorem{corollary}[theorem]{Corollary}
\newtheorem{lemma}[theorem]{Lemma}

\newtheorem{conjecture}[theorem]{Conjecture}

\newtheorem{prop}[theorem]{Proposition}
\theoremstyle{definition}

\newcommand{\ex}{\mathrm{ex}}

\newcommand{\floor}[1]{\left\lfloor #1 \right\rfloor}
\newcommand{\ceil}[1]{\left\lceil #1 \right\rceil}

\newcommand{\arrowabove}[1]{%
\begin{tikzpicture}[baseline]
    \node[anchor=base, inner sep=0pt, minimum size=0pt] (P) at (0,0) {$#1$} ;
    \useasboundingbox (current bounding box.south west) rectangle (current bounding box.north east) ;
    
    \path (P.north) ++(0ex,0.15ex) node[inner sep=0pt, minimum size=0pt, anchor=south] {$\scaleto{\rightarrow}{2pt}$} ;
\end{tikzpicture}%
}

\newcommand{\ordgraph}[3]{%
\begin{tikzpicture}[baseline]
    \node[anchor=base, inner sep=0pt, minimum size=0pt] (P) at (0,0) {${\arrowabove{#1}}^{#2}_{#3}$} ;
\end{tikzpicture}%
}

\newcommand{\Prs}{\ordgraph{P}{(r)}{s}}
\newcommand{\Krs}{\ordgraph{K}{(r)}{s}}
\newcommand{\Krn}{\ordgraph{K}{(r)}{n}}
\newcommand{\Crs}{\ordgraph{C}{(r)}{s}}

\newcommand{\UP}[2]{P^{#1}_{#2}}
\newcommand{\UPrs}{\UP{(r)}{s}}

\newcommand{\comp}[1]{\overline{#1}}
\newcommand{\st}{\colon\,}
\renewcommand{\emptyset}{\varnothing}

\newcommand{\oex}{\arrowabove{\ex}}
\newcommand{\otau}{\arrowabove{\tau}}
\newcommand{\onu}{\arrowabove{\nu}}
\newcommand{\intchi}{\chi_{i}}

\newcommand{\polylog}{\mathrm{polylog}}

\newcommand{\G}{\mathcal{G}}

\newcommand{\comment}[1]{}

\title{Tur\'an Numbers of Ordered Tight Hyperpaths}
\author{John P. Bright, Kevin G. Milans, Jackson Porter}

\begin{document}

\maketitle

\begin{abstract}
An \emph{ordered hypergraph} is a hypergraph $G$ whose vertex set $V(G)$ is linearly ordered.  We find the Tur\'an numbers for the $r$-uniform $s$-vertex tight path $\Prs$ (with vertices in the natural order) exactly when $r\le s < 2r$ and $n$ is even; our results imply $\oex(n,\Prs)=(1-\frac{1}{2^{s-r}} + o(1))\binom{n}{r}$ when $r\le s<2r$.  When $r\ge 2s$, the asymptotics of $\oex(n,\Prs)$ remain open.  For $r=3$, we give a construction of an $r$-uniform $n$-vertex hypergraph not containing $\Prs$ which we conjecture to be asymptotically extremal. 
\end{abstract}

\section{Introduction}
The \emph{Tur\'an number} of an $r$-uniform hypergraph $H$, denoted $\ex(n,H)$, is the maximum number of edges in an $r$-uniform $n$-vertex graph $G$ that does not contain $H$ as a subgraph.  Bounding Tur\'an numbers is a classical problem in extremal graph theory.  The best known general bounds on the Tur\'an numbers of the $r$-uniform $s$-vertex complete hypergraph $K^{(r)}_s$ are $(1-(\frac{r-1}{s-1})^{r-1} - o(1))\binom{n}{r} \le \ex(n,K^{(r)}_s) \le (1-\binom{s-1}{r-1}^{-1}+o(1))\binom{n}{r}$, with lower bound due to Sidorenko~\cite{SidoTuranLB} and upper bound due to de Caen~\cite{deCaenTuranUB}.

An \emph{ordered} hypergraph is a hypergraph $G$ whose vertices are linearly ordered.  For an ordered hypergraph $G$, the \emph{underlying hypergraph} is the ordinary hypergraph obtained from $G$ by discarding the order on $V(G)$.  For vertices $u$ and $v$ in an ordered hypergraph $G$, we write $u <_G v$, or $u < v$ when $G$ is clear from context, if $u$ appears before $v$ in the ordering of $V(G)$.  If $G$ and $H$ are ordered hypergraphs, then $G$ is a \emph{subgraph} of $H$, denoted $G\subseteq H$, if there is an injection $f\st V(G)\to V(H)$ such that $u<_G v$ if and only if $f(u) <_H f(v)$ and $e\in E(G)$ implies $f(e)\in E(H)$, where $f(e)=\{f(v)\st v\in e\}$.  When $H$ is an $r$-uniform ordered hypergraph, we use $\oex(n,H)$ to denote the analogous \emph{ordered Tur\'an number}, so that $\oex(n,H)$ is the maximum number of edges in an $r$-uniform $n$-vertex ordered hypergraph not containing $H$ as a subgraph.  

For graphs, ordered Tur\'an numbers behave somewhat analogously to ordinary Tur\'an numbers.  The \emph{interval chromatic number} of an ordered graph $G$, denoted $\intchi(G)$, is the minimum $k$ such that $V(G)$ can be partitioned into $k$ intervals, each of which is an independent set.  Although computing the chromatic number of an ordinary graph is NP-hard, an easy greedy algorithm computes $\intchi(G)$ for an ordered graph $G$.  Pach and Tardos~\cite{Pach2006} obtained an ordered analogue of the Erd\H{o}s--Stone Theorem, showing that for each ordered graph $H$, we have $\oex(n,H) = (1-\frac{1}{\intchi(H)-1} + o(1))\binom{n}{2}$.  Like the Erd\H{o}s--Stone Theorem, this establishes the Tur\'an numbers asymptotically for each ordered graph $G$ with $\intchi(G)>2$.  It is therefore natural to focus on ordered graphs $G$ with $\intchi(G)=2$ and ordered hypergraphs.  

A graph $G$ is a \emph{forest} if $G$ has no cycles.  Using classical Tur\'an Theory, it is straightforward to show that $\oex(n,G) \ge \Omega(n^{1+\varepsilon})$ for some positive $\varepsilon$ unless $G$ is an ordered forest with $\intchi(G)=2$.  Pach and Tardos~\cite{Pach2006} conjectured that if $G$ is an ordered forest with $\intchi(G)=2$, then $\oex(n,G) \le n(\log n)^{O(1)}$.  Kor\'andi, Tardos, Tomon, and Weidert~\cite{KorandiForest} made progress on the conjecture by proving that $\oex(n,G)\le n^{1+o(1)}$ when $G$ is an ordered forest with $\intchi(G)=2$.  For a family of ordered graphs $\G$, we define $\oex(n,\G)$ to be the maximum number of edges in an $n$-vertex ordered graph that contains no member of $\G$ as a subgraph.  A \emph{bordered cycle} is an ordered graph $G$ whose underlying graph is a cycle, whose ordering has intervals $X$ and $Y$ with $X<Y$ such that each edge in $G$ has an endpoint in $X$ and an endpoint in $Y$ (implying $\intchi(G)\le 2$), and contains the edge joining $\min X$ and $\max Y$ and the edge joining $\max X$ and $\min Y$.  Gy\H{o}ri, Kor\'andi, Methuku, Tomon, Tompkins, and Vizer~\cite{borderedcycles} proved that $\oex(n,\G_k) = \Theta(n^{1+1/k})$, where $\G_k$ is the family of bordered cycles on at most $2k$ vertices.

The $r$-uniform $s$-vertex \emph{natural path}, denoted $\Prs$, has vertex set $\{v_1,\ldots,v_s\}$ in the natural order $v_1 < \cdots < v_s$ with $E(\Prs)$ consisting of all intervals of size $r$.  The underlying hypergraph of $\Prs$ is the well-known \emph{tight path} $\UPrs$.  The \emph{length} of a path is the number of edges in the path, and so both $\Prs$ and $\UPrs$ have length $s-r+1$.  A special case of a conjecture by Kalai~\cite{KalaiConj} states that for $n\ge r\ge 2$ and $s\ge r$, we have $\ex(n,\UPrs) \le \frac{s-r}{r}\binom{n}{r-1}$, which remains open.  F\"uredi, Jiang, Kostochka, Mubayi, and Verstra\"ete~\cite{ExactResultTightPath} proved that $\ex(n,\UP{(3)}{6}) = \binom{n}{2}$ for $n\ge 5$.  In a later paper~\cite{TightHyperPaths}, the same authors proved that if $s\ge r$, then $\ex(n,\UPrs) \le \frac{s-r}{2}\binom{n}{r-1}$ when $r$ is even and $\ex(n,\UPrs) \le \frac{1}{2}(s-r+1 + \floor{\frac{s-r}{r}})\binom{n}{r-1}$ when $r$ is odd.  

\newcommand{\Qrs}{\ordgraph{Q}{(r)}{s}}
Few results on Tur\'an numbers of ordered hypergraphs are known.  In classical Tur\'an theory, an $r$-uniform hypergraph $G$ satisfies $\ex(n,G) = o(n^r)$ if and only if $G$ is \emph{$r$-partite}, meaning that there is a partition of $V(G)$ into $r$ parts such that each edge in $G$ has one vertex in each part.  The analogous statement holds for ordered hypergraphs: an $r$-uniform ordered hypergraph $G$ satisfies $\oex(n,G)=o(n^r)$ if and only if $G$ is \emph{$r$-interval-partite}, meaning that $V(G)$ can be partitioned into $r$ intervals such that each edge in $G$ has one vertex in each interval.  

For $s>r$, the natural paths $\Prs$ are not $r$-interval-partite, and so $\oex(n,\Prs) \ge \Omega(n^r)$.  The vertices of a tight path can be arranged in a different order to give an ordered $r$-interval-partite hypergraph.  The $r$-uniform $s$-vertex \emph{crossing path}, denoted $\Qrs$, is a tight path whose vertices are ordered as follows.  Arrange the $s$ vertices in an a grid with $r$ rows  $R_1,\ldots,R_r$ and $\ceil{s/r}$ columns such that any empty cells form a suffix of the last column.  Let $t=\ceil{s/r}$, and let $C_1, \ldots, C_t$ be the columns of the grid.  The ordering on the vertices of $\Qrs$ satisfies $R_1 < \ldots < R_r$, where the vertices in each $R_i$ are ordered from $C_1$ to $C_t$ (or $C_{t-1}$ if $R_i$ has no vertex in row $C_t$).  The edges of $\Qrs$ are the intervals of size $r$ in the alternative vertex ordering such that $C_1 < \cdots < C_t$, where the vertices in each $C_j$ are ordered from $R_1$ to $R_r$ (or, in the case of $C_t$, from $R_1$ to the last occupied row).  Since each edge in $\Qrs$ has one vertex in each row and $R_1 < \ldots < R_r$, it follows that $\Qrs$ is $r$-interval-partite.  F\"uredi, Jiang, Kostochka, Mubayi, and Verstra\"ete~\cite{ExtremalProbsOrderedHypergraphs} proved that $\oex(n,Q^{(r)}_s)$ equals $\binom{n}{r}-\binom{n-(s-r)}{r}$ when $r\le s\le 2r$ and is $\Theta(n^{r-1}\log n)$ when $s>2r$.  Since $\binom{n}{r} - \binom{n-(s-r)}{r} = (1+o(1))r(s-r)n^{r-1}$, it follows that always $\oex(n,Q^{(r)}_s) = O(n^{r-1}\log(n))$.  A hypergraph $F$ is a \emph{forest} if the edges of $F$ can be ordered as $e_1,\ldots,e_m$ such that for each $i$, the edge $e_i$ is the union of a subset of an earlier edge $e_j$ with $j<i$ and vertices that are not contained in any edge in $\{e_1,\ldots,e_{i-1}\}$.  In classical Tur\'an theory, we have $\ex(n,F)\le O(n^{r-1})$ for each $r$-uniform forest $F$.  
\comment{
Indeed, if $G_0$ is an $n$-vertex graph with $\Omega(n^{r-1})$ edges, we may find a copy of $F$ by reserving a set $E_1$ of at most $\binom{n}{r-1}$ edges with the property that every set of $(r-1)$ vertices contained in an edge in $G_0$ is contained in an edge in $E_1$.  Setting $G_1 = G_0 - E_1$ and iterating this process a constant $t$ number of times to produce sets $E_1,\ldots E_t$ and subgraphs $G_0,\ldots,G_t$ we may find a copy of $F$ by using any edge in $G_t$ to play the role of $e_1$ in $F$ and using edge sets $E_t,\ldots,E_1$ to find candidates for the other edges in $F$.}  Generalizing the Pach-Tardos conjecture, F\"uredi et al~\cite{ExtremalProbsOrderedHypergraphs} conjectured that $\oex(n,F) = O(n^{r-1}\cdot \polylog(n))$ when $F$ is an $r$-uniform $r$-interval-partite forest.  

We are interested in $\oex(n,\Prs)$.  In \Cref{sec:exact}, we obtain $\oex(n,\Prs)$ exactly when $s\le 2r-1$ and $n$ is even, implying that $\oex(n,\Prs)=(1-\frac{1}{2^{s-r}} + o(1))\binom{n}{r}$ when $r\le s\le 2r-1$.  When $s \ge 2r$, determining the asymptotics of $\oex(n,\Prs)$ remains open.  When $r$ divides $s$, our fractional results in  \Cref{sec:frac} imply $\oex(n,\Prs)\le (1-(\frac{r}{s})^r+o(1))\binom{n}{r}$.  When $r-1$ divides $s-1$, partitioning an interval of $n$ vertices into $(s-1)/(r-1)$ parts of equal size and removing edges with all vertices in a single part shows that $\oex(n,\Prs) \ge (1-(\frac{r-1}{s-1})^{r-1} - o(1))\binom{n}{r}$, matching the Sidorenko lower bound on $\ex(n,\Krs)$ even though $\oex(n,\Prs) \le \oex(n,\Krs) = \ex(n,K^{(r)}_s)$.  When $r$ and $s$ do not have convenient divisibility relationships, obtaining bounds on $\oex(n,\Prs)$ may involve additional subtleties.  Sidorenko's lower bound on $\ex(n,\Krs)$ holds for general $r$ and $s$; in fact, the argument shows that $\oex(n,\Crs) \ge (1-(\frac{r-1}{s-1})^{r-1} - o(1))\binom{n}{r}$, where $\Crs$ is the $r$-uniform $s$-vertex ordered tight cycle, with vertex set $\{v_0,\ldots,v_{s-1}\}$ in the natural order and edge set $\{e_0,\ldots,e_{s-1}\}$, where $e_j = \{v_j,\ldots,v_{j+r-1}\}$ (subscript arithmetic modulo $s$).

We study the case $r=3$ in \Cref{sec:r3}.  When $s$ is odd, we have that $r-1$ divides $s-1$ and the same construction as above gives $\oex(n,\Prs)\ge (1-\frac{4}{(s-1)^2}-o(1))\binom{n}{3}$.  When $s$ is even, we give a construction that improves the lower bound to $\oex(n,\Prs)\ge (1-\frac{4}{s(s-2)}-o(1))\binom{n}{3}$.  We conjecture that these bounds are asymptotically sharp.  An ordered hypergraph $G$ is \emph{monotone} if, for each edge $\{u,v,w\}\in E(G)$ with $u<v<w$, we have $\ell(uv)\le \ell(vw)$, where $\ell(xy)$ is the length of a longest ordered tight path whose last two vertices are $x$ followed by $y$ (see \Cref{sec:r3} for an equivalent, but perhaps more natural, formulation).  As some partial evidence for this conjecture, we show that if $s$ is odd and $G$ is a monotone $n$-vertex ordered hypergraph not containing $\Prs$, then $|E(G)| \le (1-\frac{4}{(s-1)^2}+o(1))\binom{n}{3}$.  The first unresolved case is that of $P^{(3)}_6$, with best known bounds $(\frac{5}{6} - o(1))\binom{n}{3} \le \oex(n,P^{(3)}_6) \le (\frac{7}{8} + o(1))\binom{n}{3}$.

\section{Exact Results for Short Paths}\label{sec:exact}

In this section, our aim is to establish $\oex(n,\Prs)$ exactly when $r\le s\le 2r-1$ and $n$ is even.  If $G \subseteq \Krn$ and $G$ does not contain $H$, then each copy of $H$ in $\Krn$ has some edge in $\comp{G}$. An \emph{$H$-transversal} in $\Krn$ is a graph $G'\subseteq \Krn$ such that every copy of $H$ in $\Krn$ has at least one edge in $G'$.  The \emph{transversal number} of $H$, denoted $\otau(n,H)$, is the minimum number of edges in an $H$-transversal.  Note that $\oex(n,H) + \otau(n,H) = |E(\Krn)| = \binom{n}{r}$.

We use $[n]$ for the vertex set of $\Krn$.  For vertex sets $A$ and $B$ in an ordered graph $G$, we write $A<B$ if $a<b$ for all $a\in A$ and $b\in B$. The \emph{reflection} of a vertex $u$ in $\Krn$ is the vertex $n+1-u$.  An \emph{interval partition} of $\Krn$ is a list of disjoint vertex subsets $(X_1,\ldots,X_k)$ whose union is $[n]$ such that each $X_i$ is an interval in $[n]$ and $X_i < X_j$ when $i<j$.  A set of vertices $S$ is \emph{$m$-left-biased} if $\Krn$ has an interval partition $(X,Y,Z)$ such that $|X| = |Z|$, $|X\cap S| = m$, and $|Z\cap S| = 0$.  Similarly, $S$ is \emph{$m$-right-biased} if $\Krn$ has an interval partition $(X,Y,Z)$ such that $|X| = |Z|$, $|X\cap S| = 0$, and $|Z\cap S| = m$.  We say that $S$ is \emph{$m$-biased} if $S$ is $m$-left-biased or $m$-right-biased.  Let $h(n,t,m)$ be the number of $t$-sets that are $m$-left-biased in $\Krn$.  Note that for even $n$, summing the $m$-left-biased $t$-sets whose $m$th vertex is at index $k$ shows that $h(n,t,m) = \sum_{k=m}^{n/2} \binom{k-1}{m-1}\binom{n-2k}{t-m}$ when $1\le m\le t$.

\newcommand{\LPrs}{\ordgraph{LP}{(r)}{s}}
Our next theorem gives a lower bound on $\oex(n,\Prs)$ by constructing an $\Prs$-transversal when $r\le s\le 2r-1$.  In fact, we construct a $\LPrs$-transversal, where $\LPrs$ is the \emph{loose path} obtained from $\Prs$ by removing all but the first and last edges.  Note that the two edges in $\LPrs$ intersect when $r<s\le 2r-1$.

\begin{theorem}\label{thm:upper}
Let $n$ be even, let $r\le s\le 2r-1$, and let $m=|E(\Prs)|=s-r+1$.  We have $\otau(n,\Prs) \le \otau(n,\LPrs) \le 2h(n,r,m) + h(n,r-1,m)$.
\end{theorem}
\begin{proof}
The condition $r\le s\le 2r-1$ translates to $1\le m\le r$.  Let $G$ be the subgraph of $\Krn$ such that $E(G) = E_1 \cup E_2$, where $E_1$ is the family of $r$-sets that are $m$-biased and $E_2$ is the family of $r$-sets whose $m$th and last vertices are reflections of one another.  Since $m>0$, the $m$-biased $r$-sets are the disjoint union of the $m$-left-biased $r$-sets and the $m$-right-biased $r$-sets, both of which have size $h(n,r,m)$, and so $|E_1| = 2h(n,r,m)$.  Also, when $m<r$, removing the last vertex from an $r$-set whose $m$th and last points are reflections of one another gives an $(r-1)$-set that is $m$-left-biased and conversely, and so $|E_2| = h(n,r-1,m)$. (Note that when $r=m$, we have $h(n,r-1,m) = 0$.)

It remains to show that every copy of $\LPrs$ in $\Krn$ has an edge in $G$.  Let $Q$ be such a copy, and let $(X,Y,Z)$ be the interval partition of $\Krn$ that minimizes $|X|$ subject to $|X|=|Z|$ and $\max\{|X\cap V(Q)|, |Z\cap V(Q)|\} \ge m$.  Such a partition exists, or else $Q$ has at most $m-1$ vertices in both the left and right halves of $\Krn$, which would imply $s=|V(Q)| \le 2(m-1) = 2(s-r)$, contradicting $s\le 2r-1$.  Let $u=|X|$.  

Suppose first that $|X\cap V(Q)| = |Z\cap V(Q)| = m$.  In this case, it must be that both $u$ and its reflection $n+1-u$ are vertices in $Q$.  Note that $s = r + (m-1)$.  Deleting the last $m-1$ vertices in $Q$ gives the first edge $e\in E(Q)$ whose $m$th vertex is $u$ and whose last vertex is $n+1-u$, implying that $e\in E_2\subseteq E(G)$.  

Otherwise, one of $\{|X\cap V(Q)|,|Z\cap V(Q)|\}$ equals $m$ and the other is at most $m-1$.  We show that $Q$ has an edge in $E_1$.  Suppose that $|X\cap V(Q)| = m$ and $|Z\cap V(Q)| \le m-1$.  Let $e$ be the first edge in $Q$ (which is obtained by deleting the last $m-1$ vertices of $Q$).  Since $s \ge m+(m-1)$, none of the deleted vertices are in $X$.  It follows that $|X\cap e| = m$ and $|Z\cap e| = 0$.  So $e$ is $m$-left-biased, and therefore $e\in E_1\subseteq E(G)$.  If instead $|X\cap V(Q)| \le m-1$ and $|Z\cap V(Q)| = m$, then a similar argument shows that the last edge $e$ in $Q$ is $m$-right-biased, also implying $e\in E_1\subseteq E(G)$.
\end{proof}

Our next theorem obtains a large family of edge-disjoint copies of $\Prs$.  For an $r$-uniform ordered hypergraph $H$, the \emph{$H$-packing number}, denoted $\onu(n,H)$ is the maximum size of an edge-disjoint family of copies of $H$ in $\Krn$.  Clearly, $\onu(n,H) \le \otau(n,H)$.

\begin{theorem}\label{thm:lower}
Let $n$ be even, let $r\le s\le 2r-1$, and let $m=|E(\Prs)|=s-r+1$.  We have $\onu(n,\Prs) \ge 2h(n,r,m) + h(n,r-1,m)$.
\end{theorem}
\begin{proof}
As in \Cref{thm:upper}, let $G$ be the subgraph of $\Krn$ with $E(G)=E_1\cup E_2$, where $E_1$ is the family of $r$-sets that are $m$-biased and $E_2$ is the family of $r$-sets whose $m$th vertex and last vertex are reflections of one another.  For each $e\in E(G)$, we construct a copy $Q_e$ of $\Prs$ such that the family of paths $\{Q_e\st e\in E(G)\}$ is edge-disjoint.

Let $e\in E(G)$.  We construct $Q_e$ as follows.  Note that every edge in $G$ is $(m-1)$-biased.  The \emph{canonical interval partition} of $e$ is the interval partition $(X,Y,Z)$ of $\Krn$ that maximizes $|X|$ subject to $|X| = |Z|$, $\max\{|X\cap e|,|Z\cap e|\} = m-1$, and $\min\{|X\cap e|, |Z\cap e|\} = 0$.  Note that by maximality of $|X|$, it follows that either $\min Y$ or $\max Y$ is a vertex in $e$, and so the canonical interval partition of $e$ is also the interval partition $(X,Y,Z)$ of $\Krn$ that minimizes $|Y|$ subject to $|X|=|Z|$ and $|Y\cap e| = r -(m-1) = 2r-s$.  The \emph{translation} of a set of vertices $S\subseteq V(\Krn)$ by a constant $c$, denoted $S+c$, is the set $\{u+c\st u\in S\}$.  If $|X\cap e| = m-1$, then we take $V(Q_e)$ to be the union of $e$ and the translation $(X\cap e) + (n-|Z|)$.  Otherwise if $|Z\cap e| = m-1$, then we take $V(Q_e)$ to be the union of $e$ and the translation $(Z\cap e) - (n-|X|)$.  In both cases, $|V(Q)| = r + m - 1 = s$, and so $Q_e$ is a copy of $\Prs$.  The \emph{core} of $Q_e$ is the set of $2r-s$ vertices in $V(Q_e)$ that belong to every edge in $Q_e$.  Since $2r-s\ge 1$, the core of $Q_e$ is non-empty.  Also, since $2r-s=s-2(m-1)$ and $|X\cap V(Q_e)| = |Z\cap V(Q_e)| = m-1$, it follows that the core of $Q_e$ equals $Y\cap V(Q_e)$.  Because the core of $Q_e$ consists of the $2r-s$ vertices in $Q$ that are closest to the center of $V(\Krn)$, given any edge in $Q_e$ we can identify the core of $Q_e$.  Moreover, since the canonical interval partition of $e$ is the partition $(X,Y,Z)$ given by minimizing $|Y|$ subject to $|X|=|Z|$ and $|Y\cap e|=2r-s$, each edge in $Q_e$ also determines the canonical interval partition of $e$.

We show that given an edge $f$ in some path $Q_e$ in the family $\{Q_e\st e\in E(G)\}$, we can determine the edge $e\in E(G)$ that generates $Q_e$.  It follows that the family is edge-disjoint.  Let $f$ be an edge in one of the paths in our collection, and recall that $f$ determines the canonical interval partition $(X,Y,Z)$ of the generating edge $e$.  It follows that $V(Q_e)$ is the union of $f$ and the translations $(X\cap f) + (n-|Z|)$ and $(Z\cap f) - (n-|X|)$.  Note that the edge $e\in E(G)$ that generates $Q_e$ must be the first or last edge in $Q_e$.  We identify $e$ as follows.  If $\min Y$ and $\max Y$ are in $V(Q_e)$, then $e\in E_2$ and $e$ is the first edge in $Q_e$.  If $\min Y$ is in $V(Q_e)$ but $\max Y$ is not, then $e\in E_1$ and $e$ is also the first edge in $Q_e$.  Otherwise, $\max Y$ is in $V(Q_e)$ and $\min Y$ is not, in which case $e\in E_1$ and $e$ is the last edge in $Q_e$.  
\end{proof}

The theorems give exact results on $\onu(n,\Prs)$, $\otau(n,\Prs)$, and $\ex(n,\Prs)$ when $r\le s\le 2r-1$.

\begin{corollary}
Let $n$ be even, let $r\le s\le 2r-1$, and let $m=|E(\Prs)| = s-r+1$.  We have that each parameter in $\{\onu(n,\Prs),\onu(n,\LPrs),\otau(n,\Prs),\otau(n,\LPrs)\}$ equals $2h(n,r,m) + h(n,r-1,m)$, and $2h(n,r,m) + h(n,r-1,m) = \frac{1}{2^{m-1}}\binom{n}{r} + O(n^{r-1})$.  Therefore $\oex(n,\LPrs) = \oex(n,\Prs) = \binom{n}{r} - 2h(n,r,m) - h(n,r-1,m) = \left(1-\frac{1}{2^{m-1}}\right)\binom{n}{r} + O(n^{r-1})$.
\end{corollary}
\begin{proof}
Clearly $\onu(n,\Prs) \le \onu(n,\LPrs), \otau(n,\Prs) \le \otau(n,\LPrs)$.  By \Cref{thm:lower} and \Cref{thm:upper} respectively, we have $\onu(n,\Prs) \ge 2h(n,r,m) + h(n,r-1,m)$ and $\otau(n,\LPrs) \le 2h(n,r,m) + h(n,r-1,m)$ and the exact results on $\onu$, $\otau$, and $\oex$ follow.

For the asymptotic results, it suffices to show that $h(n,r,m) = (1/2^m)\binom{n}{r} + O(n^{r-1})$.  Recall that $h(n,r,m)$ is the number of $m$-left-biased $r$-sets.  An $r$-set $R$ is \emph{degenerate} if $R$ contains a vertex $u$ and its reflection $n+1-u$, and $R$ is \emph{typical} if it is not degenerate.  Let $A$ be the family of degenerate $r$-sets, and let $B$ be the family of typical $r$-sets.  Note that $|A| \le (r-1)\binom{n}{r-1} = O(n^{r-1})$, since choosing $r-1$ vertices and a vertex to reflect determines a degenerate $r$-set.  Note that $|B| = 2^r\binom{n/2}{r}$, since each $r$-set in $B$ is generated by choosing $r$ of the reflection pairs, and then selecting a vertex from each chosen pair.  The $m$-left-biased sets in $B$ are obtained by choosing the left vertex from the $m$ outermost reflection pairs, and so $2^{r-m}\binom{n/2}{r}$ sets in $B$ are $m$-left-biased.  Let $C$ be the family of $m$-left-biased $r$-sets.  We compute $h(n,r,m) = |C| = |C\cap A| + |C\cap B| = O(n^{r-1}) + 2^{r-m}\binom{n/2}{r} = O(n^{r-1}) + \frac{1}{2^m}\binom{n}{r}$.
\end{proof}

\section{Fractional Variants}\label{sec:frac}
The transversal and packing numbers from \Cref{sec:exact} have  fractional variants.  For an ordered hypergraph $H$, a  \emph{fractional transversal} of the copies of $H$ in $\Krn$ is a function $w$ that assigns non-negative weights to each edge in $\Krn$ such that $\sum_{e\in E(H')} w(e) \ge 1$ for each copy $H'$ of $H$ in $\Krn$.  The \emph{fractional transversal number} of $H$, denoted $\otau^*(n, H)$, is the infimum, over all fractional transversals $w$, of the sum of $w(e)$ over all edges $e\in E(\Krn)$.  Standard compactness arguments show that the infimum is always achieved, and so we may replace infimum by minimum in the definition.  Also, if $G'\subseteq \Krn$ and $G'$ is an $H$-transversal, then the weight function $w$ with $w(e)=1$ for $e\in E(G')$ and $w(e)=0$ for $e\not\in E(G')$ is a fractional transversal, and therefore $\otau^*(n, H) \le \otau(n,H)$.  

Each fractional transversal of copies of $H$ in $\Krn$ is a feasible solution to the linear program with variables $\{w(e)\st e\in E(\Krn)\}$ with objective to minimize $\sum_e w(e)$ subject to $\sum_{e\in E(H')} w(e)\ge 1$ for each copy $H'$ of $H$ in $\Krn$.  The dual linear program has variables $\{w(H')\st\mbox{$H'$ is a copy of $H$ in $\Krn$}\}$ with objective to maximize $\sum_{H'} w(H')$ subject to the constraints that, for each edge $e$ in $\Krn$, the sum of $w(H')$ over all copies $H'$ of $H$ in $\Krn$ that contain $e$ is at most $1$.  A feasible solution $w$ to the dual linear program is called a \emph{fractional $H$-packing}, and the \emph{fractional $H$-packing number}, denoted $\onu^*(n,H)$ is the value of this linear program.  Since both the LP and its dual are clearly feasible, it follows from theory of linear programming that $\otau^*(n,H) = \onu^*(n,H)$.  As before, standard compactness arguments show that a fractional $H$-packing with total weight $\onu^*(n,H)$ exists, and it is clear that $\onu(n,H)\le \onu^*(n,H)$.  Hence 
\[ \onu(n,H) \le \onu^*(n,H) = \otau^*(n,H) \le \otau(n,H).\]
In this section, we show that $\onu(n,\Prs)$, $\onu^*(n,\Prs)$, and $\otau^*(n,\Prs)$ are all asymptotically $\left(\left(\frac{r}{s}\right)^r + o(1)\right)\binom{n}{r}$ when $r$ divides $s$.

\begin{prop}\label{prop:flower}
If $r$ divides $s$, then $\onu(n,\Prs) \ge \left(\left(\frac{r}{s}\right)^r + o(1)\right)\binom{n}{r}$.
\end{prop}
\begin{proof}
We give a $\Prs$-packing of the required size.  Let $k=s/r$, and without loss of generality assume that $k$ divides $n$.  Let $(X_1,\ldots,X_k)$ be an interval partition of $\Krn$ into parts of equal size.  For each edge $e\in E(\Krn)$ with $e\subseteq X_1$, let $P_e$ be the $s$-vertex path with vertex set $\bigcup_{j=0}^{k-1} (e + j|X_1|)$.  Note that given any edge $e'\in E(P_e)$, it is easy to recover $e$, and it follows that $\{P_e\st e\subseteq X_1\}$ is an edge-disjoint collection of copies of $\Prs$.  Therefore $\onu(n,\Prs) \ge \binom{|X_1|}{r} = \binom{n/k}{r} = \left(\frac{1}{k^r} + o(1)\right)\binom{n}{r}$.
\end{proof}

\begin{prop}\label{prop:fupper}
If $r$ divides $s$, then $\otau^*(n,\Prs)\le \left(\left(\frac{r}{s}\right)^r + o(1)\right)\binom{n}{r}$.
\end{prop}
\begin{proof}
We give a fractional $\Prs$-transversal.  Let $k=s/r$.  We may assume without loss of generality that $k$ divides $n$.  Let $(X_1,\ldots,X_k)$ be an interval partition of $\Krn$ into parts of equal size.  Let $w$ be the weight function with $w(e)=r/s$ if $e$ is contained in a part in $(X_1,\ldots,X_k)$ and $w(e)=0$ otherwise.  Let $P$ be a copy of $\Prs$ in $\Krn$.  Note that at most $r-1$ vertices in $V(P)\cap X_i$ begin an edge with weight zero, and $P$ has at most $k(r-1)$ such vertices.  Therefore at least $s-k(r-1)$ vertices in $P$ begin an edge with positive weight.  So the edges of $P$ have total weight at least $(r/s)(s - k(r-1))$, and this equals $1$.  It follows that $w$ is a fractional $\Prs$-transversal, and so $\otau^*(n,\Prs) \le \frac{r}{s}\cdot k\cdot\binom{|X_1|}{r} = \binom{n/k}{r} = \left(\frac{1}{k^r} + o(1)\right)\binom{n}{r}$.
\end{proof}

If $r$ divides $s$ and the integer $s/r$ also divides $n$, then \Cref{prop:flower} and \Cref{prop:fupper} imply $\onu(n,\Prs) = \otau^*(n,\Prs)=\onu^*(n,\Prs)=\binom{n/k}{r}$, where $k=s/r$. 

\begin{theorem}
If $r$ divides $s$, then $\onu(n,\Prs)$, $\onu^*(n,\Prs)$, and $\otau^*(n,\Prs)$ are all asymptotically equal to $\left(\left(\frac{r}{s}\right)^r + o(1)\right)\binom{n}{r}$.
\end{theorem}

The smallest path to which our argument in \Cref{sec:exact} does not apply is $P^{(r)}_{2r}$, when $s=2r$.  In this case, our fractional results imply $\otau^*(n,P^{(r)}_{2r}) = \left(\frac{1}{2^r} + o(1)\right)\binom{n}{r}$, but, at least in the case $r=3$, we believe that $\otau(n,P^{(r)}_{2r}) \gg \otau^*(n,P^{(r)}_{2r})$.  In particular, $\otau^*(n,P^{(3)}_6) = \left(\frac{1}{8} + o(1)\right)\binom{n}{r}$ but we conjecture $\otau(n,P^{(3)}_6) = \left(\frac{1}{6} + o(1)\right)\binom{n}{r}$.

\section{The case $r=3$}\label{sec:r3}
\newcommand{\Pthrees}[0]{P^{(3)}_{s}}
\newcommand{\Ktwon}{\ordgraph{K}{(2)}{n}}
When $r=3$, the ordered Tur\'an numbers are equivalent to an edge-labeling problem on the ordered complete graph $\Ktwon$.  A \emph{$k$-edge-labeling} of $\Ktwon$ assigns to each pair $uv$ a label in a linearly ordered set $S$ with $|S|=k$.  Let $\phi$ be a $k$-edge-labeling of $\Ktwon$.  A triple of vertices $\{u,v,w\}$ with $u<v<w$ is \emph{good} if $\phi(uv) < \phi(vw)$.  For convenience, we write $uvw$ for the triple $\{u,v,w\}$.  A triple is \emph{bad} if it is not good.  Let $f(n,k)$ be the maximum, over all $k$-edge-labelings $\phi$ of $\Ktwon$, of the number of good triples.  A $k$-edge-labeling $\phi$ of $\Ktwon$ is \emph{optimal} if it has $f(n,k)$ good triples.

\begin{prop}\label{prop:labelequiv}
For $s\ge 3$, we have $\oex(n,\Pthrees) = f(n,s-2)$.
\end{prop}
\begin{proof}
First, we show $\oex(n,\Pthrees)\ge f(n,s-2)$.  Let $\phi$ be an $(s-2)$-edge-labeling of $\Ktwon$ with $f(n,s-2)$ good triples.  Let $G$ be the ordered $3$-uniform hypergraph with vertex set $V(\Ktwon)$ such that for $u<v<w$, we have $uvw\in E(G)$ if and only if $\phi(uv)<\phi(vw)$.  We claim that $G$ does not contain $\Pthrees$.  Indeed, if $v_1\cdots v_s$ is a copy of $\Pthrees$ in $G$, then $\phi(v_{i-1}v_i) < \phi(v_iv_{i+1})$ for $1<i<s$ by definition of $G$.  It follows that $\phi$ uses $s-1$ distinct labels on the consecutive pairs of $v_1\cdots v_s$, contradicting that $\phi$ is a $(s-2)$-edge-labeling.  It follows that $\oex(n,\Pthrees)\ge |E(G)| = f(n,s-2)$.

Next, we show $\oex(n,\Pthrees)\le f(n,s-2)$.  Let $G$ be a $3$-uniform ordered hypergraph not containing $\Pthrees$, and let $\phi$ be the edge-labeling on $V(G)$ by setting $\phi(uv)$, where $u<v$, equal to the length of a longest tight ordered path in $G$ that ends in $uv$.  Clearly, if $u<v<w$ and $uvw\in E(G)$, then we have $\phi(uv)<\phi(vw)$ since the edge $uvw$ can be used to extend a longest ordered path ending in $uv$ to a longest path of larger length ending in $vw$.  Therefore $\phi$ has at least $|E(G)|$ good triples.  Note that $\phi$ assigns each pair $uv$ a value in the set $\{0,\ldots,s-3\}$, since every ordered tight path in $G$ has at most $s-1$ vertices and at most $s-3$ edges.  Since $\phi$ assigns each edge a value in $\{0,\ldots,s-3\}$, it follows that $\phi$ is an $(s-2)$-edge-labeling, and therefore $f(n,s-2) \ge |E(G)| = \oex(n,\Pthrees)$.
\end{proof}

Next, we give lower bound constructions for $f(n,k)$ which we conjecture to be asymptotically optimal.  The construction is easiest to describe when $k$ is odd.  A labeling $\phi$ is \emph{monotone} if, for all $u<v<w$, we have $\phi(uv)\le\phi(vw)$.

\begin{prop}\label{prop:labelodd}
If $k$ is odd, then $f(n,k)\ge \left(1-\frac{4}{(k+1)^2} + o(1)\right)\binom{n}{3}$.
\end{prop}
\begin{proof}
Let $t=(k+1)/2$.  Let $(X_1,\ldots,X_t)$ be an interval partition of $\Ktwon$ into $t$ parts whose sizes differ by at most $1$.  For $u<v$ with $u\in X_i$ and $v\in X_j$, we set $\phi(uv) = i+j$.  Clearly, the range of $\phi$ is contained in $\{2,\ldots,2t\}$, and so $\phi$ is a  $(2t-1)$-edge-labeling.  Note that $2t-1=k$.

Since $\phi$ is a monotone labeling, the only triples $uvw$ with $u<v<w$ that are not good are those with $\phi(uv)=\phi(vw)$.  Each such triple is contained in a part $X_i$ for some $i$.  It follows that the number of triples that are not good is asymptotically equal to $t\binom{n/t}{3}$, which is asymptotically equal to $\frac{1}{t^2}\binom{n}{3}$.  The proposition follows.
\end{proof}

\begin{corollary}
If $s\ge 3$ and $s$ is odd, then $\oex(n,\Pthrees) \ge \left(1-\frac{4}{(s-1)^2}+o(1)\right)\binom{n}{3}$.
\end{corollary}

Applying \Cref{prop:labelequiv} to the construction in \Cref{prop:labelodd} gives a graph $G$ whose complement is the union of $t$ complete graphs on $t$ disjoint intervals of nearly equal size.  The construction for even $k$ is more subtle.

We first give our construction in terms of a general interval partition $(X_1,\ldots,X_k)$ of $[n]$ into $k$ parts.  Later, we specify the sizes of the parts.  For a pair $uv$ with $u<v$, $u\in X_i$, and $v\in X_j$, we define $\phi(uv)$ as follows.  If $i=j$, then $\phi(uv)=i$.  If $j-i \ge 2$, then we set $\phi(uv)$ so that $i<\phi(uv)<j$.  Otherwise, $j=i+1$.  The \emph{fractional index} of $u$ in $X_i$ is $(u+1-\min X_i)/|X_i|$.  Note that the fractional index of $u$ is a real number in $(0,1]$.  Let $\lambda_u$ and $\lambda_v$ be the fractional indices of $u$ in $X_i$ and $v$ in $X_{i+1}$, respectively.  We set $\phi(uv)=i$ if $\lambda_u + \lambda_v \le 1$ and $\phi(uv)=i+1$ if $\lambda_u + \lambda_v > 1$.  

\begin{lemma}\label{lem:badcount}
Let $a,b,c$ be constants.  If $|X_{i-1}| = (a+o(1))n$, $|X_i| = (b+o(1))n$, and $|X_{i+1}|=(c+o(1))n$, then the number of bad triples $uvw$ with $v\in X_i$ is $\left[\left(a+b\right)b\left(b+c\right) + o(1)\right]\binom{n}{3}$.
\end{lemma}
\begin{proof}
\comment{
Let $X_{i-1} = \{u_1,\ldots,u_r\}$, $X_i = \{v_1,\ldots,v_s\}$, and let $X_{i+1} = \{w_1,\ldots,w_t\}$.  
}
Let $v\in X_i$ and let $\lambda_v$ be the fractional index of $v$ in $X_i$.  If $uvw$ is a bad triple with $u<v<w$, then $u\in X_{i-1} \cup X_i$, $w\in X_i\cup X_{i+1}$, and $\phi(uv)=\phi(vw)=i$.  For $u\in X_i$, we require only that $u$ precede $v$, and there are $\lambda_v|X_i| - 1$ such vertices in $X_i$.  For $u\in X_{i-1}$, the condition $\phi(uv) = i$ is equivalent to $\lambda_u + \lambda_v > 1$, where $\lambda_u$ is the fractional index of $u\in X_{i-1}$.  The number of $u\in X_{i-1}$ with $\lambda_u + \lambda_v > 1$ is $|X_{i-1}| - \floor{(1-\lambda_v)|X_{i-1}|}$ or $\ceil{\lambda_v|X_{i-1}|}$.  It follows that the number of choices for $u$ in a bad triple $uvw$ with $u<v<w$ equals $\lambda_v(|X_{i-1}| + |X_i|) + O(1)$.

Similarly, the number of $w\in X_i$ that follow $v$ is $(1-\lambda_v)|X_i|$ and the number of $w\in X_{i+1}$ with $\lambda_v + \lambda_w \le 1$ is $\floor{(1-\lambda_v)|X_{i+1}|}$.  It follows that the number of choices for $w$ in a bad triple $uvw$ with $u<v<w$ equals $(1-\lambda_v)(|X_i| + |X_{i+1}|) + O(1)$.

Multiplying the number of choices for $u$ and the number of choices for $w$ gives a total of $\lambda_v(1-\lambda_v)(|X_{i-1}| + |X_i|)(|X_i| + |X_{i+1}|) + O(n)$ bad triples $uvw$ with $u<v<w$.  Suppose that $X_i = \{v_1,\ldots,v_t\}$.  Summing over all $v\in X_i$, the total number of bad triples is $O(n^2) + (|X_{i-1}| + |X_i|)(|X_i| + |X_{i+1}|)\frac{1}{t^2} \sum_{j=1}^t j(t-j)$, or $O(n^2) + (|X_{i-1}| + |X_i|)(|X_i| + |X_{i+1}|)\frac{1}{t^2} \binom{t+1}{3}$.  Recalling that $t=|X_i| = (b+o(1))n$, the number of bad triples $uvw$ with $u<v<w$ simplifies to $O(n^2) + (a+b)(b+c)b(\frac{1}{6} + o(1))n^3$ and the lemma follows.
\end{proof}

\begin{lemma}\label{lem:labeleven}
If $k$ is even, then $f(n,k)\ge \left(1-\frac{4}{k(k+2)} - o(1)\right)\binom{n}{3}$.
\end{lemma}
\begin{proof}
Suppose $k$ is even, and let $t=k/2$.  Let $(Y_1,\ldots,Y_t)$ and $(Z_1,\ldots,Z_{t+1})$ be interval partitions of $V(\Ktwon)$ into parts of nearly equal size, and let $(X_1,\ldots,X_k)$ be their common refinement.  (Note that $(Y_1,\ldots,Y_t)$ is the partition used in our construction with $k-1$ labels, and $(Z_1,\ldots,Z_{t+1})$ is the partition used in our construction with $k+1$ labels.)  For $1\le j\le k$, we set $a_j$ equal to the limit of $|X_j|/n$ as $n\to\infty$.  It is convenient to introduce $X_0 = X_{k+1} = \emptyset$ and $a_0 = a_{k+1} = 0$.  When divided by $n$ to normalize, the boundaries of $(Y_1,\ldots,Y_t)$ are $\frac{0}{t},\frac{1}{t},\ldots,\frac{t}{t}$, and the boundaries of $(Z_1,\ldots,Z_{t+1})$ are $\frac{0}{t+1},\frac{1}{t+1},\ldots,\frac{t+1}{t+1}$.  In the common refinement, the these boundaries interleave and are thus $\frac{0}{t},\frac{1}{t+1},\frac{1}{t},\frac{2}{t+1},\frac{2}{t},\ldots,\frac{t}{t+1},\frac{t}{t}$.  It follows that for $0\le j\le t$, we have $a_{2j} = \left(\frac{j}{t} - \frac{j}{t+1}\right) = \frac{j}{t(t+1)}$, and for $0\le j\le t$, we have $a_{2j+1} = \left(\frac{j+1}{t+1}-\frac{j}{t}\right) = \frac{t-j}{t(t+1)}$.

Let $\phi$ be the labeling described above.  Since $k$ is constant, by \Cref{lem:badcount}, the number of bad triples is asymptotically $\sum_{i=1}^k \left(a_{i-1} + a_i\right)a_i\left(a_i + a_{i+1}\right)\binom{n}{3}$.  Let $A$ be this sum taken over even indices $i$, and let $B$ be the sum taken over odd indices $i$.  We compute
\begin{align*}
A &= \sum_{j=1}^t \left(a_{2j-1} + a_{2j}\right)a_{2j}\left(a_{2j} + a_{2j+1}\right)\cdot\binom{n}{3} \\
&= \binom{n}{3}\sum_{j=1}^t \frac{(t-(j-1))+j}{t(t+1)} \cdot \frac{j}{t(t+1)} \cdot \frac{j+(t-j)}{t(t+1)} \\
&= \binom{n}{3}\frac{1}{\left(t(t+1)\right)^2}\sum_{j=1}^t j\\
&= \frac{1}{2t(t+1)}\binom{n}{3}
\end{align*}
The computation for $B$ is symmetric and also gives $B=\frac{1}{2t(t+1)}\binom{n}{3}$.  Since the total number of bad triples is asymptotic to $A+B$, we have $f(n,k) \ge \left(1-\frac{1}{t(t+1)} - o(1)\right)\binom{n}{3}$.  
\end{proof}

\begin{corollary}
If $s\ge 2$ and $s$ is even, then $\oex(n,\Pthrees) \ge \left(1-\frac{4}{(s-2)s}+o(1)\right)\binom{n}{3}$.
\end{corollary}

We conjecture that these constructions are asymptotically optimal.

\begin{conjecture}
If $a=\floor{(k+1)^2/4}$, then $f(n,k) = \left(1 - \frac{1}{a} + o(1)\right)\binom{n}{3}$.  Equivalently, if $b=\floor{(s-1)^2/4}$, then $\oex(n,\Pthrees)=\left(1-\frac{1}{b}+o(1)\right)\binom{n}{3}$.
\end{conjecture}

\newcommand{\val}{\mathrm{cost}}
Our goal in the remainder of this section is to show that if $k$ is odd and $\phi$ is monotone, then $\phi$ has at least $\left(\frac{4}{(k+1)^2} - o(1)\right)\binom{n}{3}$ bad triples.  This shows that the construction in \Cref{prop:labelodd} is asymptotically optimal within the class of monotone labelings.  For an edge-labeling $\phi$, let $\val(\phi)$ be the number of bad triples.

We believe that there exists an optimal monotone labeling of $\Ktwon$.  It is at least true that when $k\ge 2$, an optimal $[k]$-edge-labeling $\phi$ of $\Ktwon$ does not contain a triple $uvw$ with $u<v<w$, $\phi(uv)=k$ and $\phi(vw)=1$.  Indeed, if there is such a triple, then we may fix $v$ and assume that $u$ is the minimum vertex such that $\phi(uv)=k$ and $w$ is the maximum vertex such that $\phi(vw)=1$. Modify $\phi$ to obtain $\phi'$ by setting $\phi'(uv)=1$ and $\phi'(vw)=k$.  We have $\val(\phi')=\val(\phi)+a-b$, where $a$ is the number of triples that are good in $\phi$ but bad in $\phi'$, and $b$ is the number of triples that are bad in $\phi$ but good in $\phi'$.  We show that $b>a$, contradicting that $\phi$ is optimal.  By maximality of $w$, for each $w'$ with $w'>w$, the triple $uvw'$ contributes to $b$.  Similarly by minimality of $u$, every triple $u'vw$ with $u'<u$ also contributes to $b$.  Also, the triple $uvw$ itself contributes to $b$.  Therefore $b\ge (n-w) + (u-1) + 1$.  Every triple contributing to $a$ contains $uv$ (and therefore has the form $u'uv$ for some $u'<u$), or $vw$ (and therefore has the form $vww'$ for some $w'>w$).  It follows that $a\le (u-1)+(n-w) < (n-w)+(u-1)+1\le b$.

For a monotone $[k]$-edge-labeling $\phi$ of $\Ktwon$ define $\Phi_L$ and $\Phi_R$ as follows.  We set $\Phi_L(1) = 0$ and $\Phi_L(v) = \max \{\phi(uv)\st u < v\}$ for $v>1$.  We set $\Phi_R(n)=k+1$ and $\Phi_R(v) = \min\{\phi(vw) \st w > v\}$ for $v<n$. Note that by monotonicity of $\phi$, if $u < v$, then $\Phi_L(u) \le \Phi_R(u) \le \phi(uv) \le \Phi_L(v) \le \Phi_R(v)$.

Using $\Phi_L$ and $\Phi_R$, we construct two interval partitions of $V(\Ktwon)$.  For $0\le i\le k+1$, let $X_i = \{v\in V(\Ktwon)\st \Phi_L(v)=i\}$ and let $\hat X_i = \{v\in V(\Ktwon)\st \Phi_R(v)=i\}$.  Note that $X_0=\{1\}$, $\hat X_0 = \emptyset$, $X_{k+1} = \emptyset$, and $\hat X_{k+1} = \{n\}$.  Since $\Phi_L$ and $\Phi_R$ inherit the monotonicity of $\phi$, both $(X_0,\ldots,X_{k+1})$ and $(\hat X_0,\ldots,\hat X_{k+1})$ are interval partitions of $V(\Ktwon)$.  Our next lemma shows that these partitions are very similar.

\begin{lemma}\label{lem:partitions}
The symmetric difference of $X_i$ and $\hat X_i$ has size at most $2$, and so $||X_i| - |\hat X_i|| \le 2$.
\end{lemma}
\begin{proof}
Let $u,v\in X_i$ with $u<v$.  Since $i = \Phi_L(u) \le \Phi_R(u)\le \phi(uv)\le \Phi_L(v)=i$ it follows that $\Phi_R(u) = i$, and so $u\in \hat X_i$.  Therefore $X_i - \hat X_i \subseteq \{\max X_i\}$.  Similarly, $\hat X_i - X_i \subseteq \{\min \hat X_i\}$.
\end{proof}

A key step in our proof is the following bound on the sizes of the parts. 

\begin{lemma}\label{lem:leftparts}
Let $\phi$ be an monotone $[k]$-edge-labeling of $\Ktwon$ that minimizes the number of bad triples, and define $\Phi_L$ and $(X_0,\ldots,X_{k+1})$ as above.  For $2\le i\le k$, we have $|X_i| + |X_{i+1}| \le |X_{i-2}| + |X_{i-1}| + 2$.
\end{lemma}
\begin{proof}
First, suppose that $X_i$ is nonempty.  Let $v=\min X_i$ and let $u$ be the least vertex such that $\phi(uv)=i$.  (Note that $u$ exists since $\Phi_L(v)=i>0$.)  If it exists, let $w$ be the least vertex in $X_{i+1}$ with $\phi(vw) = i+1$.  We say that an edge is \emph{long} if its endpoints are in distinct non-consecutive parts in $(X_0,\ldots,X_{k+1})$. 

Obtain $\phi'$ from $\phi$ by reducing by $1$ the labels on $vw$ (if $w$ exists), $uv$, all long edges $u'u$ such that $u'<u$ and $\phi(u'u) = i-1$, and all long edges $v'v$ such that $v'<v$ and $\phi(v'v)=i$.  Note that $\phi'$ is a monotone labeling.  We have $\val(\phi') = \val(\phi) + a - b$, where $a$ is the number of triples $xyz$ with $x<y<z$ which are bad in $\phi'$ and good in $\phi$, and $b$ is the number of triples $xyz$ with $x<y<z$ which are good in $\phi'$ and bad in $\phi$.  Note that if $xyz$ is a triple with $x<y<z$ and $\phi$ and $\phi'$ agree on both $xy$ and $yz$, then of course $xyz$ contributes to neither $a$ nor $b$.  Also, if $\phi$ and $\phi'$ disagree on $xy$ and $yz$, then also $xyz$ does not contribute to $a$ or $b$ since $\phi'(xy)=\phi(xy)-1$ and $\phi'(yz)=\phi(yz)-1$.  So each triple $xyz$ contributing to $a$ or $b$ has one pair where $\phi$ and $\phi'$ agree, and one pair where $\phi$ and $\phi'$ disagree.  Suppose that $xyz$ contributes to $a$.  It follows that $\phi'(xy)=\phi(xy)$ and $\phi'(yz)=\phi(yz)-1 = \phi(xy)$.  Note that $yz$ is not a long edge $u'u$, since $\phi(xu') \le \Phi_L(u') \le \Phi_L(u) - 2 \le i-3$ and $\phi(u'u) = i-1$, implying that $xu'u$ is still good in $\phi'$.  Similarly, $yz$ is not a long edge $v'v$.  So $yz\in \{uv,vw\}$.  

If $xyz=xuv$, then $\phi(xu)=\phi'(xu)=\phi'(uv)=i-1$ and hence $u \in X_{i-1}$.  It follows that $x\in X_{i-2}\cup X_{i-1}$ since $xu$ is not a long edge.  So $xyz=xuv$ implies that $x<u$ and $x\in X_{i-2}\cup X_{i-1}$.

If $xyz=xvw$, then $\phi(xv) = \phi'(xv) = \phi'(vw)=i$.  Since $u$ is the minimum vertex with $\phi(uv)=i$, we have $u\le x$.  Moreover, $xyz\ne uvw$ since $\phi$ and $\phi'$ both disagree on $uv$ and $vw$.  Hence $u<x<v$. Also, $x\in X_{i-1}$ since $\phi(xv)=\phi'(xv)$, and so $xv$ is not a long edge with label $i$.  Combining both cases, the contributions $xyz$ to $a$ arise from distinct $x\in X_{i-2}\cup X_{i-1}$, and it follows that $a\le |X_{i-2}| + |X_{i-1}|$.  

It remains to show that $b\ge |X_i|+ |X_{i+1}|-2$.  Let $z\in X_i \cup X_{i+1}$ such that $z>v$ and $z\ne w$ (if $w$ exists).  Note that $\phi(vz) = \phi'(vz)$, since $vw$ is the only edge incident to the right of $v$ where $\phi$ and $\phi'$ disagree.  If $\phi(vz)=i$, then $uvz$ contributes to $b$ since $\phi(uv)=i$ but $\phi'(uv)=i-1$. If $\phi(vz) > i$, then $i < \phi(vz) \le \Phi_L(z) \le i+1$ and so $z\in X_{i+1}$.  Hence $w$ exists and by minimality of $w$ we have that $w < z$. Since $\phi(vw) = i+1$ and also $i+1=\Phi_L(w) \le \phi(wz) \le \Phi_L(z) \le i+1$, we have also $\phi(wz)=i+1$.  But $\phi'(vw)=i$ and $\phi'(wz)=\phi(wz)=i+1$, and it follows that $vwz$ contributes to $b$.  Therefore $b\ge |X_i| + |X_{i+1}| - 2$.

Since $\phi$ minimizes the number of bad triples among monotone labelings, it follows that $\val(\phi')\ge \val(\phi)$, giving $a\ge b$.  

Finally, we consider the case that $X_i = \emptyset$.  If $X_{i+1} = \emptyset$ also, then the inequality holds trivially.  Let $v=\min\{X_{i+1}\}$.  Since $\Phi_L(v)=i+1$, there is a vertex $u$ such that $u<v$ and $\phi(uv)=i+1$.  Since $X_i=\emptyset$, it follows that $\Phi_L(u)\le i-1$, implying that $uv$ is a long edge.  Obtain $\phi'$ from $\phi$ by reducing $\phi(uv)$ by $1$, and leaving all other labels the same.  As before, we have $\val(\phi')=\val(\phi)+a-b$, where $a$ is the number of triples $xyz$ that are good in $\phi$ and bad in $\phi'$, and $b$ is the number of triples $xyz$ that are bad in $\phi$ and good in $\phi'$.  A contribution $xyz$ to $a$ requires that $yz=uv$ and $\phi'(xu) = \phi'(uv) = i$, but $\phi'(xu) = \phi(xu) \le \Phi_L(xu) \le i-1$, and so $a=0$.  Also, if $z>w$ and $z\in X_{i+1}$, then $uvz$ contributes to $b$ since $\phi'(vz)=\phi(vz)=i+1$ and $\phi'(uv) = \phi(uv) - 1 = i$.  It follows that $b\ge |X_{i+1}| - 1$ and by optimality of $\phi$, we obtain $0 = a\ge b \ge |X_{i+1}| - 1$, and the inequality follows.
\end{proof}

\begin{corollary}\label{cor:rightparts}
Let $\phi$ be a monotone $[k]$-edge-labeling of $\Ktwon$ that minimizes the number of bad triples, and define $\Phi_R$ and $(\hat X_0,\ldots,\hat X_{k+1})$ as above.  For $1\le j\le k-1$, we have $|\hat X_{j-1}| + |\hat X_j| \le |\hat X_{j+1}| + |\hat X_{j+2}|+2$.
\end{corollary}
\begin{proof}
Obtain $\phi'$ from $\phi$ by reversing the order of vertices in $\Ktwon$ and inverting the labels, so that $\phi'(uv)=(k+1) - \phi(v^*u^*)$, where $v^* = (n+1)-v$ and $u^*=(n+1)-u$.  Note that $\phi'$ is a monotone $[k]$-edge-labeling with $\val(\phi')=\val(\phi)$.  Moreover, defining $\Phi'_L$ with respect to $\phi'$ and the corresponding partition $(X'_0,\ldots,X'_{k+1})$, we have that $\Phi'_L(u) = (k+1)-\Phi_R(u^*)$ where $u^*=(n+1)-u$ and $|X'_i| = |\hat X_{(k+1)-i}|$.  Applying \Cref{lem:leftparts} to $\phi'$ with $i=(k+1)-j$ gives the result.
\end{proof}

\begin{theorem}\label{thm:monobound}
Let $\phi$ be a monotone $[k]$-edge-labeling of $\Ktwon$ that minimizes the number of bad triples.  If $k$ is odd, then $\val(\phi) \ge (1-o(1))\frac{4}{(k+1)^2}\binom{n}{3}$.
\end{theorem}
\begin{proof}
Define $\Phi_L$, $\Phi_R$,  $(X_0,\ldots,X_{k+1})$, and $(\hat X_0,\ldots,\hat X_{k+1})$ as above, and let $m=(k-1)/2$.  For $0\le \ell \le m$, let $a_\ell = |X_{2\ell}| + |X_{2\ell + 1}|$ and let $\hat a_\ell= |\hat X_{2\ell}| + |\hat X_{2\ell + 1}|$.  By \Cref{lem:leftparts} with $i=2\ell$, we have $a_{\ell-1} \ge a_\ell - 2$ for $1\le \ell \le m$.  It follows that $a_0 \ge a_1 - 2 \ge \cdots \ge a_m - 2m = a_m - (k-1)$.  With two applications of \Cref{lem:partitions}, we have $a_0 \le \hat a_0 + 4$.  By \Cref{cor:rightparts} with $j=2\ell-1$, we have $\hat a_{\ell-1} \le \hat a_\ell + 2$ for $1\le \ell \le m$, and it follows that $\hat a_0 \le \hat a_1 + 2 \le \cdots \le \hat a_m + 2m = \hat a_m + (k-1)$.  Chaining the inequalities gives $-(k-1) + a_m \le \cdots \le a_0 \le \hat a_0 + 4 \le \cdots \le \hat a_m + k+3 \le a_m + k+7$.  It follows that $-2\ell + a_\ell$ is in the interval $[a_m - (k-1), a_m + (k+7)]$, and so $a_\ell$ is in the interval $[a_m - (k-1), a_m + (k+7)+2\ell]$.  Since $2\ell \le k-1$, we have $|a_\ell - a_m| \le 2k+6$ and so each pair in $\{a_0,\ldots,a_m\}$ differs by at most $4k+12$.  Since $\sum_{\ell=0}^m a_\ell = n$, it follows that $a_\ell = n/(m+1) + O(k)$.

Similarly, for $0\le \ell \le m$, let $b_\ell = |X_{2\ell + 1}| + |X_{2\ell + 2}|$ and $\hat b_\ell = |\hat X_{2\ell + 1}| + |\hat X_{2\ell + 2}|$.  By \Cref{lem:leftparts} with $i=2\ell+1$, we  have $b_{\ell - 1} \ge b_\ell - 2$ for $1\le \ell\le m$.  It follows that $b_0 \ge b_1 - 2 \ge \cdots \ge b_m - (k-1)$.  Similarly, $b_0 \le \hat b_0 + 4$.  Also, \Cref{cor:rightparts} with $j=2\ell$ gives $\hat b_{\ell-1} \le \hat b_{\ell} + 2$ for $1\le \ell\le m$.  Therefore $\hat b_0 \le \hat b_1 + 2 \le \cdots \le \hat b_m + (k-1)$.  Combining the inequalities gives $-(k-1) + b_m \le \cdots \le b_0 \le \hat b_0 + 4 \le \hat b_m + (k+3) \le b_m + (k+7)$.  As above, $b_\ell$ and $b_{\ell'}$ differ by at most $4k+12$.  Again, $\sum_{\ell=0}^m b_\ell = n$, and so $b_\ell = n/(m+1) + O(k)$.  

Let $0\le \ell \le m$.  We claim that $X_{2\ell}$ and $X_{2\ell + 2}$ differ in size by at most $O(k)$.  Indeed, both $a_\ell$ and $b_\ell$ equal $n/(m+1) + O(k)$, and so $|a_\ell - b_\ell| \le O(k)$.  Since $a_\ell = |X_{2\ell}| + |X_{2\ell + 1}|$ and $b_\ell = |X_{2\ell + 1}| + |X_{2\ell + 2}|$, we have that $a_\ell - b_\ell = |X_{2\ell}| - |X_{2\ell + 2}|$ and the claim follows.  Therefore each pair of parts in $\{X_0,X_2,X_4,\ldots,X_{k+1}\}$ differs in size by at most $O(k^2)$.  Since $|X_{k+1}| = 0$, it follows that each part with even index has size at most $O(k^2)$.  Since $a_l = |X_{2\ell}|+|X_{2\ell+1}|=n/(m+1)+O(k)$, it follows that each part with odd index has size $n/(m+1)-O(k^2)$.  Since each triple $uvw$ with $u,v,w\in X_i$ satisfies $\phi(uv)=\phi(vw)=i$, the number of bad triples in $\phi$ is at least $\sum_{\ell = 0}^m \binom{|X_{2\ell + 1}|}{3}$, and $\sum_{\ell = 0}^m \binom{|X_{2\ell + 1}|}{3} \ge (m+1)\binom{(n/(m+1)) - O(k^2)}{3} = (1-o(1))\frac{1}{(m+1)^2}\binom{n}{3}$.
\end{proof}

Let $k$ be even, let $\phi$ be an optimal monotone $[k]$-edge-labeling of $\Ktwon$, and define the parts $(X_0,\ldots,X_{k+1})$ as above.  Similar arguments as in \Cref{thm:monobound} can be used to obtain the sizes of the parts asymptotically, and these match the sizes of the corresponding parts in our construction in \Cref{lem:labeleven}.  However, bounding the number of bad triples in $\phi$ for even $k$ is more complicated since consecutive parts $X_i$ and $X_{i+1}$ are both linear in $n$, making the triples with two vertices in one of $\{X_i,X_{i+1}\}$ and one in the other significant.

\bibliographystyle{plainnat}
\bibliography{citations}

\end{document}